\documentclass[12pt,leqno]{article}

\usepackage{amsmath,amssymb,amsthm}
\usepackage[all]{xy}
\usepackage{lscape}

\makeatletter

\swapnumbers

\newtheorem{theorem}{Theorem}[section]
\newtheorem{proposition}[theorem]{Proposition}
\newtheorem{lemma}[theorem]{Lemma}
\newtheorem{corollary}[theorem]{Corollary}

\theoremstyle{definition}

\newtheorem{example}[theorem]{Example}

\theoremstyle{remark}

\theoremstyle{remark}
\newtheorem{remark}[theorem]{Remark}

\def\({{\rm (}}
\def\){{\rm )}}

\let\Mathrm\operator@font
\let\Cal\mathcal
\let\Bbb\mathbb
\newcommand{\fm}{\ensuremath{\mathfrak m}}

\def\standop#1{\mathop{\Mathrm #1}\nolimits}
\def\difstop#1#2{\expandafter\def\csname #1\endcsname{\standop{#2}}}
\def\defstop#1{\difstop{#1}{#1}}

\defstop{AB}
\defstop{ann}
\defstop{Ass}

\defstop{codim}
\defstop{Coh}
\defstop{Coker}
\defstop{Cone}

\defstop{depth}

\defstop{EM}
\defstop{End}
\defstop{ev}
\defstop{Ext}

\defstop{Flat}
\defstop{Func}

\def\GL{{\sl{GL}}}

\difstop{height}{ht}
\defstop{Hom}
\let\hst\heartsuit

\def\id{\mathord{\Mathrm{id}}}
\def\Id{\mathord{\Mathrm{Id}}}
\difstop{Image}{Im}

\defstop{Ker}

\defstop{Lch}
\defstop{length}
\defstop{Lin}
\defstop{Lqc}
\defstop{lqc}

\defstop{Mat}
\defstop{Min}
\defstop{Mod}
\defstop{Mor}

\defstop{Nerve}

\def\op{^{\standop{op}}}

\defstop{PA}
\defstop{PM}
\defstop{Proj}

\defstop{Qch}
\defstop{qch}

\defstop{rank}

\defstop{res}

\def\Sch{\underline{\Mathrm Sch}}

\defstop{Spec}
\defstop{supp}
\defstop{Supp}
\defstop{Sym}

\difstop{tdeg}{trans.deg}

\defstop{Tor}
\difstop{trace}{tr}

\defstop{Zar}

\def\I{\Cal I}
\def\O{\Cal O}
\def\C{\Cal C}
\def\M{\Cal M}
\def\N{\Cal N}

\def\P{\Cal P}
\def\Q{\Cal Q}
\def\E{\Cal E}
\def\fH{\mathfrak{H}}
\def\fG{\mathfrak{G}}
\def\fD{\mathfrak{D}}
\def\fd{\mathfrak{d}}
\def\fM{\mathfrak{M}}
\def\fm{\mathfrak{m}}

\let\indlim\varinjlim
\let\projlim\varprojlim

\def\sHom{\mathop{{}^*\Mathrm{Hom}}\nolimits}
\def\uHom{\mathop{\mbox{\underline{$\Mathrm Hom$}}}\nolimits}
\def\uExt{\mathop{\mbox{\underline{$\Mathrm Ext$}}}\nolimits}
\def\uGamma{\mathop{\mbox{\underline{$\Gamma$}}}\nolimits}
\def\uH{\mathop{\mbox{\underline{$H$}}}\nolimits}

\def\sdarrow#1{\downarrow\hbox to 0pt{\scriptsize$#1$\hss}}
\def\suarrow#1{\uparrow\hbox to 0pt{\scriptsize$#1$\hss}}
\def\ssearrow#1{\searrow\hbox to 0pt{\scriptsize$#1$\hss}}

\def\grad#1{\sqrt[G]{#1}}

\def\section{\@startsection{section}{1}{\z@ }%
{-3.5ex plus -1ex minus -.2ex}{2.3ex plus .2ex}{\bf }}

\long\def\refname{\par\kern -3ex
\begin{center}\rm R\sc{eferences}\end{center}\par\kern 
-2ex}

\def\@seccntformat#1{\csname the#1\endcsname.\quad}

\def\@@@sect#1#2#3#4#5#6[#7]#8{%
   \ifnum #2>\c@secnumdepth 
      \def \@svsec {}\else \refstepcounter {#1}%
      \def\@svsec{}
   \fi 
   \@tempskipa #5\relax 
   \ifdim \@tempskipa >\z@ 
     \begingroup #6\relax \@hangfrom {\hskip #3\relax 
     \@svsec}{\interlinepenalty \@M #8\par }\endgroup 
     \csname #1mark\endcsname {#7}
   \else 
   \def \@svsechd {#6\hskip #3\@svsec #8\csname #1mark\endcsname {#7}}
   \fi \@xsect {#5}}

\def\@@@startsection#1#2#3#4#5#6{%
 \if@noskipsec \leavevmode \fi \par \@tempskipa #4\relax \@afterindenttrue 
 \ifdim \@tempskipa <\z@ \@tempskipa -\@tempskipa \@afterindentfalse 
 \fi \if@nobreak \everypar {}\else \addpenalty {\@secpenalty }\addvspace 
  {\@tempskipa }\fi \@ifstar {\@ssect {#3}{#4}{#5}{#6}}{\@dblarg 
  {\@@@sect {#1}{#2}{#3}{#4}{#5}{#6}}}}

\def\theparagraph{\thesection.\arabic{paragraph}}
\def\aparagraph{\@@@startsection{paragraph}{2}{\z@ }%
              {1.75ex plus .2ex minus .15ex}{-1em}{\bf(\theparagraph) } }
\def\paragraph{\@@@startsection{paragraph}{2}{\z@ }%
              {1.75ex plus .2ex minus .15ex}{-1em}{}{\bf(\theparagraph)} }

\let\c@theorem\c@paragraph

\title{Equivariant Matlis and the local duality}

\author{M{\sc itsuyasu} H{\sc ashimoto} and M{\sc asahiro} O{\sc htani}}
\date{\normalsize
Graduate School of Mathematics, Nagoya University\\
Chikusa-ku,  Nagoya 464--8602 JAPAN\\
{\small \tt hasimoto@math.nagoya-u.ac.jp~~~
m05011w@math.nagoya-u.ac.jp}}

\begin{document}

\maketitle

\begin{abstract}
Generalizing the known results on graded rings and modules,
we formulate and prove the equivariant version of the local duality on
schemes with a group action.
We also prove an equivariant analogue of Matlis duality.
\end{abstract}

\section{Introduction}

This paper is a continuation of \cite{HO}, and study equivariant 
local cohomology.
In this paper, utilizing an equivariant dualizing complex, we define the
$G$-sheaf of matlis, an equivariant analogue of the injective hull of the
residue field of a local ring.
Using this, we formulate and prove Matlis and the local duality under
equivariant settings.

Let $R$ be a Gorenstein local ring, $T=R[x_1,\ldots,x_s]$ be the graded
polynomial ring with $r_i:=\deg x_i$ positive, $I$ a homogeneous ideal of 
height $h$, and $A:=T/I$.
Assume that $A$ is Cohen--Macaulay of dimension $d$.
Set $\omega_T:=T(-r)$, where $r=\sum_i r_i$, and $(-r)$ denotes the 
shift of degree.
Set $\omega_A:=\Ext_T^h(A,\omega_T)$.
For a graded $A$-module $M$, set 
$M^\vee:=\bigoplus_{i\in\Bbb Z}M_{-i}^\dagger$, where $(?)^\dagger=
\Hom_R(?,E_R)$, where $E_R$ is the injective hull of the residue field of
$R$.
Note that $M^\vee$ is a graded $A$-module again.
Note also that $\sHom_A(M,A^\vee)\cong M^\vee$ (see for the notation
$\sHom$, \cite[page~33]{BH}).

For a finite graded $A$-module $M$, 
we have an isomorphism of {\em graded} $A$-modules
\[
H^i_{\fM}(M)\cong \Ext^{d-i}_A(M,\omega_A)^\vee,
\]
cf.\ \cite[Theorem~3.6.19]{BH}, 
see also Corollary~\ref{graded-local-duality.thm}.

The main purpose of this paper is to generalize this graded version of 
local duality to more general equivariant local duality.
Note that a graded module over a $\Bbb Z$-graded ring is nothing but
an equivariant module under the action of $\Bbb G_m=\GL_1$,
see \cite[(II.1.2)]{Hashimoto}.
On the way, we prove some basic properties on equivariant local cohomology.

In this introduction, let $S$ be a noetherian scheme, 
$G$ a flat $S$-group scheme of finite type, and $X$ a noetherian $G$-scheme.
In order to establish an analogy of the local duality on $X$, we need to 
define an equivariant analogue of a local ring or a local scheme.
This is done in \cite{HO}, and it is a $G$-local $G$-scheme.
So let $X$ be a $G$-local $G$-scheme.
That is to say, $X$ has a unique minimal nonempty $G$-stable closed
$G$-subscheme, say $Y$.
Next, we need to have an equivariant analogue of local cohomology.
This is the main subject of \cite{HO}.
Finally, we need to have an analogy of the Matlis duality.
In other words, we need to have an analogue of the injective hull of the 
residue field of a local ring.
The authors do not know how to define it quite generally.
However, if $X$ has a $G$-equivariant dualizing complex (see for the 
definition, \cite[chapter~31]{ETI}) $\Bbb I_X$, then we can define it
as the unique nonzero cohomology group of $R\uGamma_Y(\Bbb I_X)$.
We call this sheaf the $G$-sheaf of Matlis.
Thus we can formulate the equivariant local duality.
The proof depends on the isomorphism $\fH$, see below.

Many ideas used in this paper have already appeared in 
the theory of graded rings
\cite{GW1}, \cite{GW2}, \cite{BH}, \cite{Kamoi}.
If $H$ is a finitely generated abelian group, then letting 
$G= \Spec \Bbb ZH$, where $\Bbb ZH$ is the group algebra of $H$ over $\Bbb Z$,
an $H$-graded algebra is nothing but a $G$-algebra, and for a $G$-algebra
$A$, a graded $A$-module is nothing but a $(G,A)$-module.
However, we need to point out that for a general $G$ and a $G$-local
$G$-algebra $(A,\fM)$ with the $G$-dualizing complex $\Bbb I$, 
the global section of the $G$-sheaf of Matlis $E_A$ is not necessarily
injective as a $(G,A)$-module, see Example~\ref{non-injective.thm}.
In particular, $E_A$ is not the injective hull of $A/\fM$ in the category
of $(G,A)$-modules.


Using the $G$-sheaf of Matlis, 
we can prove a weak version of the Matlis duality, too.
It is a duality from the category of coherent $(G,\O_X)$-modules of finite
length to itself, see Theorem~4.17.
Note that a better Matlis duality exists over a complete local ring.
It is a duality from the category of noetherian modules to the category of
artinian modules (\cite[Theorem~3.2.13]{BH}).
The authors do not know a good analogue of a complete local ring, and
thus cannot give an equivariant Matlis duality between noetherian 
quasi-coherent $(G,\O_X)$-modules and artinian modules in general.
However, there is an example of graded case of that kind of duality, 
see Remark~\ref{equivariant-Matlis-graded.thm}.

Section~2 is preliminaries.
We give some basic properties of the duality map in a closed category.
We also give some sufficient conditions to guarantee that injective objects
in the category $\Qch(G,X)$ is acyclic with respect to some cohomological
functors.
We also prove a generalization of the flat base change 
(\cite[Theorem~6.10]{HO}), see Lemma~\ref{flat-base-change.thm}.
We also describe the local cohomology over a diagram of schemes
using the inductive limit of $\uExt$ groups, as in the single-scheme case.

Section~3 treats the map $\fH$.
For a small category $I$, an $I\op$-diagram of schemes $X$, an open 
subdiagram of schemes $U$ of $X$, and an open subdiagram of schemes 
$V$ of $U$, there is a natural map
\[
\fH: \uGamma_{U,V}\uHom_{\O_X}(\M,\N)\rightarrow 
\uHom_{\O_X}(\M,\uGamma_{U,V}\N)
\]
for $\M,\N\in\Mod(X)$.
There is an obvious derived version of it, and $\fH$ is often an isomorphism
(see Lemma~\ref{fH-isom.thm} and  Theorem~\ref{fH-isom-derived.thm}).
This is the key to the proof of the equivariant version of the
local duality.
In order to establish the existence and some basic properties of $\fH$,
we need to prove various commutativity of diagrams.
To do this, we utilize the basics on closed categories as in 
\cite[chapter~1]{ETI}.

In section~4, we formulate and prove the equivariant analogues of
Matlis and the local duality.
We start with Matijevic--Roberts type theorem for $G$-local $G$-schemes,
and prove an equivariant version of Nakayama's lemma, which is well-known
for affine case.

In section~5, we give an example of the graded case.
Note that in some cases, Matlis duality can be in more general form than
the version described in section~4, see 
Remark~\ref{equivariant-Matlis-graded.thm}.

\section{Preliminaries}

\paragraph
We use the notation and terminology of \cite{ETI}, \cite{HO}, and \cite{HM}
freely.

\paragraph\label{fD.par}
Let $X$ be a symmetric monoidal closed category (see
\cite[(3.5.1)]{Lipman}), and $b,d\in X$.
Then we denote the composite map
\[
b \xrightarrow{\trace}[[b,d],b\otimes [b,d]]\xrightarrow{\gamma}
[[b,d],[b,d]\otimes b]\xrightarrow{\ev}[[b,d],d]
\]
by $\fD$, and we call it the {\em duality map},
where $\trace$, $\gamma$, and $\ev$ denote 
the trace map \cite[(1.30)]{ETI}, the twisting (symmetry) isomorphism 
\cite[(1.28)]{ETI}, and the evaluation map \cite[(1.30)]{ETI}, respectively.

\begin{lemma} $\fD$ is natural on $b$.
Namely, for a morphism $\phi:b\rightarrow b'$, the diagram
\[
\xymatrix{
b \ar[r]^\phi \ar[d]^{\fD} & b' \ar[d]^{\fD} \\
[[b,d],d] \ar[r]^\phi & [[b',d],d]
}
\]
is commutative.
\end{lemma}

\proof Consider the diagram
\[
\xymatrix{
 & b \ar[ld]_{\phi} \ar[d]^{\trace} \ar[r]^-{\trace}
 \ar@{}[dr]|{\text{(a)}} \ar@{}[ddl]|{\text{(c)}} &
[[b,d],b\otimes [b,d]] \ar[d]^\phi
    \ar@{}[dr]|{\text{(b)}} 
    \ar `r[rdd] `[dd]^{\ev\gamma} [dd] & \\
b' \ar[d]^{\trace} &
[[b',d],b\otimes [b',d]] \ar[ld]_{\phi} \ar[r]^\phi 
\ar@{}[d]|{\text{(d)}} &
[[b',d],b\otimes[b,d]] \ar[ld]^{\ev\gamma} & \\
[[b',d],b'\otimes [b',d]] \ar[r]^-{\ev\gamma} & 
[[b',d],d] & [[b,d],d] \ar[l]_\phi & 
}.
\]
(a) and (d) are commutative by \cite[Lemma~1.32]{ETI}.
The commutativity of (b) and (c) are trivial.
\qed

\begin{lemma}
For a morphism $\psi:d \rightarrow d'$, the diagram
\[
\xymatrix{
b \ar[r]^{\fD} \ar[d]^{\fD} &
[[b,d],d] \ar[d]^\psi \\
[[b,d'],d'] \ar[r]^\psi &
[[b,d],d']
}
\]
is commutative.
\end{lemma}

\proof Consider the diagram
\[
\xymatrix{
b \ar[r]^-{\trace} \ar[d]^{\trace} \ar@{}[dr]|{\text{(a)}} &
[[b,d],b\otimes [b,d]] \ar[d]^\psi \ar[r]^-{\ev\gamma}
\ar@{}[dr]|{\text{(b)}} &
[[b,d],d] \ar[dd]^\psi \\
[[b,d'],b\otimes [b,d']]\ar[r]^\psi \ar[d]^{\ev\gamma}
\ar@{}[dr]|{\text{(c)}} &
[[b,d],b\otimes [b,d']] \ar[dr]^{\ev\gamma} & \\
[[b,d'],d'] \ar[rr]^{\psi} & & [[b,d],d']
}.
\]
(a) is commutative by \cite[Lemma~1.32]{ETI}.
(b) and (c) are obviously commutative.
Hence the whole diagram is commutative.
\qed

\begin{lemma}\label{fD-H.thm}
Let $f:X\rightarrow Y$ be a symmetric monoidal functor 
{\rm\cite[(3.4.2)]{Lipman}}
between symmetric monoidal closed categories.
For $b,d\in X$, the diagram
\[
\xymatrix{
fb \ar[r]^-{\fD} \ar[d]^{\fD} &
f[[b,d],d] \ar[d]^H \\
[[fb,fd],fd] \ar[r]^H &
[f[b,d],fd]
}
\]
is commutative.
\end{lemma}

\proof Consider the diagram
\[
\xymatrix{
 & [[fb,fd],fb\otimes [fb,fd]] \ar[d]^H 
    \ar `l[lddd] `[ddd]^{\ev\gamma} [ddd] \ar@{}[dr]|{\text{(a)}} &
fb \ar[d]^\trace \ar[l]_-\trace 
\ar `r[ddr]^\trace[ddr] & \\
\ar@{}[dr]|{\text{(b)}} &
[f[b,d],fb\otimes[fb,fd]] \ar[d]^{\ev\gamma}
\ar@{}[dr]|{\text{(d)}} &
[f[b,d],fb\otimes f[b,d]] \ar[d]^m \ar@{}[r]|{\text{(c)}} \ar[l]_H & \\
 & [f[b,d],fd] &
[f[b,d],f(b\otimes[b,d])] \ar[l]_-{\ev\gamma} \ar@{}[dr]|{\text{(e)}} &
f[[b,d],b\otimes[b,d]] \ar[l]_-H \ar`d[dl][dl]_{\ev\gamma} \\
 & [[fb,fd],fd] \ar[u]_H &
f[[b,d],d] \ar[ul]_H &
}.
\]

(a) is commutative by \cite[(1.32)]{ETI}.
The commutativity of (b) is trivial.
(c) is \cite[(1.37)]{ETI} and is commutative.
(d) is \cite[(1.36)]{ETI} and is commutative.
(e) is commutative by the naturality of $H$.
\qed

\paragraph A symmetric monoidal functor $f:X\rightarrow Y$ is said to be
{\em Lipman} if it has a left adjoint $g:Y\rightarrow X$ such that
the natural maps $\Delta:g(b\otimes d)\rightarrow gb\otimes
gd$ and $C:g\O_Y\rightarrow \O_X$ are isomorphisms,
see \cite[(1.48)]{ETI}.
We also say that $(f,g)$ is a Lipman adjoint pair in this case.

By Lemma~\ref{fD-H.thm}, we have:

\begin{lemma}\label{fD-P.thm}
Let $f:X\rightarrow Y$ and $g:Y\rightarrow X$ be a Lipman adjoint pair
where $X$ and $Y$ are closed.
Then the diagram
\[
\xymatrix{
gb' \ar[r]^-{\fD} \ar[d]^{\fD} &
g[[b',d'],d'] \ar[d]^P \\
[[gb',gd'],gd'] \ar[r]^P &
[g[b',d'],gd']
}
\]
is commutative.
\end{lemma}

\paragraph\label{premodule-fD.par}
Let $(\Bbb X,\O_{\Bbb X})$ be a ringed category.
That is, $\Bbb X$ is a small category, and $\O_{\Bbb X}$ is a presheaf of
commutative rings on $\Bbb X$.
Then for $\M,\N\in\PM(\Bbb X)$, the map
\[
\fD: \M\rightarrow \uHom_{\PM(\Bbb X)}(\uHom_{\PM(\Bbb X)}(\M,\N),\N)
\]
is described as follows.
At $x\in \Bbb X$, 
\begin{multline*}
\fD:\Gamma(x,\M)\rightarrow 
\Gamma(x,\uHom_{\PM(\Bbb X)}(\uHom_{\PM(\Bbb X)}(\M,\N),\N))\\
=\Hom_{\PM(\Bbb X)/x}(\uHom_{\PM(\Bbb X)}(\M,\N)|_x,\N|_x))
\end{multline*}
is given as follows.
For $a\in\Gamma(x,\M)$, $\fD(a):\uHom_{\PM(\Bbb X)}(\M,\N)|_x\rightarrow
\N|_x$ is the map such that for $\phi:y\rightarrow x$, 
$\fD(a)_\phi:\Hom_{\PM(\Bbb X/y)}(\M|_y,\N|_y)\rightarrow \Gamma(y,\N)$
is given by $\fD(a)_\phi(h)=h(a)$.
This is proved easily using \cite[(2.42)]{ETI} and \cite[(2.41)]{ETI}.

\paragraph Let $(\Bbb X,\O_{\Bbb X})$ be a ringed site, 
and $\M,\N\in\Mod(\Bbb X)$.
Then the map
\[
\fD: \M\rightarrow \uHom_{\O_{\Bbb X}}(\uHom_{\O_{\Bbb X}}(\M,\N),\N)
\]
is exactly the same map as the one described in (\ref{premodule-fD.par}).
This follows from \cite[(2.49)]{ETI}, Lemma~\ref{fD-H.thm}, and 
(\ref{premodule-fD.par}).

\paragraph
In the rest of this paper, $S$ denotes a scheme, and $G$ an $S$-group scheme.
We write diagrams of schemes as $X$, $Y$, $Z,\ldots$ (not as
$X_\bullet$, $Y_\bullet$, $Z_\bullet,\ldots$).
Similarly, morphisms of diagrams of schemes are expressed as $f$, $g$, 
$h,\ldots$, not as $f_\bullet$, $g_\bullet$, $h_\bullet,\ldots$
This is a convention in \cite{HO}.

\begin{lemma}\label{lipman-indlim.thm}
Let $I$ be a small category, and 
$f:X\rightarrow Y$ be a concentrated \(i.e., quasi-compact quasi-separated\)
morphism of $I\op$-diagrams of schemes.
Let $(C_\alpha)$ be a pseudo-filtered inductive system of complexes of 
$\Cal O_X$-modules such that 
for each $j\in I$, one of the following holds:
\begin{description}
\item[(a)]
There exists some $n_j\in\Bbb Z$ such that for any $\alpha$, 
$\tau^{\leq n_j-1}(C_\alpha)_j$ is exact
\(see for the definition of $\tau^{\leq n_j-1}$, see 
{\rm\cite[(3.24)]{ETI}}\){\rm;}
\item[(b)] 
Each point of $X_j$ has a noetherian open neighborhood of
finite Krull dimension.
\item[(c)] For any $\alpha$, 
$C_{\alpha,j}$ has quasi-coherent cohomology
groups.
\end{description}
Set $C=\indlim C_\alpha$.
Then the canonical map
\begin{equation}\label{lipman.eq}
\indlim R^if_*C_\alpha \rightarrow R^if_* C
\end{equation}
is an isomorphism for $i\in\Bbb Z$.
If, moreover, each $C_\alpha$ is $f_*$-acyclic, then $C$ is $f_*$-acyclic.
\end{lemma}

\proof In view of \cite[Example~8.23, {\bf 2}]{ETI}, 
it is easy to see that it suffices to show that
\[
\indlim R^i(f_j)_*C_{\alpha,j}\rightarrow R^i(f_j)_*C_j
\]
is an isomorphism for each $j$, to prove that (\ref{lipman.eq}) is an
isomorphism.
This is (3.9.3.1) and (3.9.3.2) of \cite{Lipman}.

To prove the last assertion, it suffices to show that
each $C_j$ is $(f_j)_*$-acyclic.
This is \cite[(3.9.3.4)]{Lipman}.
\qed

\begin{corollary}\label{lqc-complex-acyclic.thm}
Let $f:X\rightarrow Y$ be as in {\rm Lemma~\ref{lipman-indlim.thm}}.
Let $C$ be a complex of $\O_X$-modules such that each term of $C$
is locally quasi-coherent and $f_*$-acyclic.
Then $C$ is $f_*$-acyclic.
\end{corollary}

\proof Similar to \cite[(3.9.3.5)]{Lipman}.
\qed

\begin{lemma}\label{K-inj-acyclic.thm}
Let $X$ and $Y$ be $S$-groupoid \(see for the definition, 
{\rm\cite[(12.1)]{ETI})}
and $f:X\rightarrow Y$ a 
morphism \(in the category $\P(\Delta_M,\Sch/S)$, 
see for the notation, {\rm\cite[Glossary]{ETI})}.
Assume that $f$ is cartesian, $Y$ has affine arrows, 
and assume one of the following:
\begin{description}
\item[(a)]
$X_0$ is noetherian;
\item[(b)] 
$Y_0$ and $f_0$ are quasi-compact separated.
\end{description}
Then
\begin{description}
\item[(i)] $f$ is concentrated and $X_0$ is concentrated.
\item[(ii)]
A $K$-injective complex $\Bbb I$ in $K(\Qch(X))$ is
$f_*$-acyclic.
\item[(iii)] The canonical maps
\[
F_Y\circ Rf_*^{\Qch}\cong
R(F_Y\circ f_*^{\Qch})\cong
R(f_*\circ F_X)\rightarrow Rf_*\circ F_X
\]
are all isomorphisms, where $F_Y:D(\Qch(Y))\rightarrow D(Y)$ and
$F_X:D(\Qch(X))\rightarrow D(X)$ are triangulated functors induced by
inclusions, and $f_*^{\Qch}:\Qch(X)\rightarrow\Qch(Y)$ is the restriction
of $f_*:\Mod(X)\rightarrow\Mod(Y)$, see {\rm\cite[Lemma~7.14]{ETI}}.
\end{description}
\end{lemma}

\proof {\bf (i)} In either case, $f_0$ is concentrated.
Since $f$ is cartesian, 
each $f_i$ ($i=0,1,2$) is obtained as a base change of $f_0$,
and hence is concentrated.
It is easy to see that $X_0$ is concentrated in either case.

{\bf (ii)} As $f$ is concentrated cartesian,
$f_*^{\Qch}$ is well-defined \cite[Lemma~7.14]{ETI}.
Since $X_0$ is concentrated and $X$ has affine arrows,
$\Qch(X)$ is Grothendieck by \cite[Lemma~12.8]{ETI}.
So $\Bbb I$ has a strictly injective resolution (that is, a $K$-injective
resolution each of whose term is injective) $\Bbb J$ 
\cite[Proposition~3.2]{Franke}.
As the mapping cone of $\Bbb I\rightarrow \Bbb J$ is null-homotopic,
replacing $\Bbb I$ by $\Bbb J$, we may assume that $\Bbb I$ is
strictly injective.
By Corollary~\ref{lqc-complex-acyclic.thm}, it suffices to show that
each term of $\Bbb I$ is $f_*$-acyclic.
So we may assume that $\Bbb I$ is a single injective object of $\Qch(X)$.
Let $\Bbb I_0\rightarrow K$ be a monomorphism with $K$ an injective object
of $\Qch(X_0)$.
This is possible, since $\Qch(X_0)$ is Grothendieck 
\cite[Corollary~11.7]{ETI}.
Note that the restriction $(?)_0:\Qch(X)\rightarrow\Qch(X_0)$ has
the right adjoint $(d_0)_*^{\Qch}\circ \Bbb A$, see \cite[Lemma~12.11]{ETI}.
As $(?)_0$ is faithful exact, the composite
\[
\Bbb I\rightarrow (d_0)_*^{\Qch}\Bbb A\Bbb I_0\rightarrow
(d_0)_*^{\Qch}\Bbb A K
\]
is a monomorphism into an injective object.
This must split, and hence we may further assume that $\Bbb I=(d_0)_*^{\Qch}
\Bbb A K$.

By restriction, it suffices to show that $R^j(f_i)_*\Bbb I_i=0$ for
$i=0,1,2$ and $j>0$.
Since $(?)_{i}\Bbb A\cong r_0(i+1)^*$ (see for the
notation, \cite[(9.1)]{ETI}) and $(d_0)$ is affine,
\begin{multline*}
R^j(f_i)_*\Bbb I_i\cong
R^j(f_i)_*d_0(i+1)_*r_0(i+1)^* K
\cong
R^j(f_i\circ d_0(i+1))_* r_0(i+1)^*K\\
=R^j(d_0(i+1)\circ f_{i+1})_*r_0(i+1)^*K
\cong d_0(i+1)_* R^j (f_{i+1})_* r_0(i+1)^*K\\
\cong
d_0(i+1)_*r_0(i+1)^*R^j (f_0)_*K
=0
\end{multline*}
for $j>0$ by \cite[Lemma~14.6, {\bf 1}]{ETI} and its proof.
This is what we wanted to prove.

{\bf (iii)} Follows immediately from {\bf (ii)}.
\qed

The following is a generalization of \cite[Theorem~6.10]{HO}.

\begin{lemma}\label{flat-base-change.thm}
Let $I$ be a small category, 
$h:X'\rightarrow X$ a flat morphism of $I\op$-diagrams of schemes.
Let $f:U\hookrightarrow X$ be an open subdiagram of schemes, and
$g:V\hookrightarrow U$ be an open subdiagram of schemes.
Assume that $f$ and $g$ are quasi-compact.
Let $f':U'\hookrightarrow X'$ and $g':V'\hookrightarrow U'$ be the
base change of $f$ and $g$, respectively.
Then $\bar\delta:h^*R\uGamma_{U,V}\rightarrow R\uGamma_{U',V'}h^*$ in
{\rm\cite[(6.1)]{HO}} is an isomorphism between functors from 
$D_{\Lqc}(X)$ to $D_{\Lqc}(X')$.
\end{lemma}

\proof As in the proof of \cite[Corollary~6.3]{HO}, we may assume that
the problem is on single schemes.
Consider the map of triangles
\[
\xymatrix{
h^*R\uGamma_{U,V} \ar[r] \ar[d]^{\bar\delta} &
h^*Rf_*f^* \ar[r] \ar[d]^{d\theta} &
h^*Rf_*Rg_*g^*f^* \ar[r] \ar[d]^{dd\theta\theta} &
h^*R\uGamma_{U,V}[1] \ar[d] \\
R\uGamma_{U',V'}h^*\ar[r] &
Rf'_*(f')^*h^* \ar[r] &
Rf'_*Rg'_*(g')^*(f')^*h^* \ar[r] &
R\uGamma_{U',V'}h^*[1]
}.
\]
By \cite[Proposition~3.9.5]{Lipman}, the vertical arrows $d\theta$ and
$dd\theta\theta$ are isomorphisms.
Hence, $\bar\delta$ is also an isomorphism.
\qed

\paragraph\label{indlim.thm}
Let $I$ be a small category, $X$ an $I\op$-diagram of
schemes, and $Y$ a cartesian closed subdiagram of schemes of $X$
defined by the quasi-coherent ideal sheaf $\I$ of $\O_X$.
Assume that $X$ is locally noetherian with flat arrows.
Then, the canonical map
\[
\Phi_Y: \indlim \uHom_{\O_X}(\O_X/\I^n,\M)\rightarrow \uGamma_Y\M
\]
is an isomorphism for $\M\in\Lqc(X)$, see \cite[(3.21)]{HO}.
By the way-out lemma \cite[Proposition~I.7.1]{Hartshorne}, we have

\begin{lemma}
Let the notation be as in {\rm(\ref{indlim.thm})}.
Then for $\Bbb F\in D_{\Lqc}^+(X)$, 
\[
\Phi_Y: R(\indlim \uHom_{\O_X}(\O_X/\I^n,?))(\Bbb F)\rightarrow R\uGamma_Y\Bbb
F
\]
is an isomorphism.
In particular, $\Phi_Y$ induces an isomorphism
\[
\indlim \uExt_{\O_X}^i(\O_X/\I^n,\Bbb F)\cong \uH^i_Y(\Bbb F).
\]
\qed
\end{lemma}

\begin{lemma}
Let $X$ be an $S$-groupoid with affine arrows.
Let $U$ be a cartesian open subdiagram of $X$,
and $V$ a cartesian open subdiagram of $Y$.
Assume that $X_0$ is noetherian.
If $\Bbb I$ is a $K$-injective complex in $K(\Qch(X))$, then 
$\Bbb I$ is $\uGamma_{U,V}$-acyclic.
\end{lemma}

\proof Using \cite[Corollary~6.7]{HO},
it suffices to show that for an injective object $K$ of $\Qch(X_0)$,
$(d_0)_*^{\Qch}\Bbb A K$ is $\uGamma_{U,V}$-acyclic,
as in the proof of Lemma~\ref{K-inj-acyclic.thm}.
Applying restrictions, 
it suffices to show that
$\uH^j_{U_i,V_i}(d_0(i+1)_*r_0(i+1)^*K)=0$ for $j>0$ and $i=0,1,2$.
By the independence \cite[Corollary~4.17]{HO} and the flat base change
Lemma~\ref{flat-base-change.thm}, 
this sheaf is $d_0(i+1)_*r_0(i+1)^*\uH^j_{U_0,V_0}K$.
Since $K$ is also injective in $\Mod(X)$ \cite[Theorem~II.7.18]{Hartshorne},
it is a flabby sheaf, and $\uH^j_{U_0,V_0}K=0$.
\qed

\paragraph
A $G$-scheme $X$ (i.e., an $S$-scheme with a left $G$-action) is said
to be {\em standard} if
$X$ is noetherian, and 
the second projection $p_2:G\times X\rightarrow X$ is
flat of finite type.

Let $X$ be a standard $G$-scheme.
We denote the category of quasi-coherent (resp.\ coherent) $(G,\O_X)$-modules
by $\Qch(G,X)$ (resp.\ $\Coh(G,X)$).
Note that the sheaf theory discussed in \cite[chapters~29--31]{ETI}
and \cite{HO},
where we assume that $G$ is flat of finite type over $S$, 
still works under our weaker assumption ($p_2$ is flat of finite type).
In particular, $\Qch(G,X)$ is a locally noetherian category, and 
$\M\in\Qch(G,X)$ is a noetherian object if and only if $\M\in\Coh(G,X)$,
see \cite[Lemma~12.8]{ETI}.

\paragraph We say that a standard $G$-scheme $X$ 
is {\em $G$-artinian} if there is no incidence
relation between $G$-prime $G$-ideals (see for the definition, 
\cite[(4.12)]{HM}) of $X$.

\begin{lemma}\label{G-artinian.thm}
If $X$ is $G$-artinian, then $X$ is a disjoint union of finitely many
$G$-artinian $G$-local $G$-schemes.
\end{lemma}

\proof Clearly, the set of all $G$-prime $G$-ideals $\Spec_G(X)$
agrees with the finite set $\Min_G(\O_X)$, 
the set of minimal $G$-primes of $0$.
Thus there are only finitely many $G$-prime $G$-ideals.
For $\P,\Q\in\Spec_G(X)$ with $\P\neq\Q$, 
$\Ass_G(\O_X/(\P+\Q))=\emptyset$, 
since there is no $G$-prime $G$-ideal containing both $\P$ and $\Q$.
Thus $\P+\Q=\O_X$.
This shows that $X=\coprod_{\P\in\Spec_G(X)}V(\P)$.
As each $V(\P)$ is clearly $G$-artinian $G$-local, we are done.
\qed

\section{The map $\fH$}

\paragraph Let $f:X\rightarrow Y$ be a symmetric monoidal functor between
symmetric monoidal closed categories, and $g:Y\rightarrow X$ its right
adjoint.
For $b\in Y$ and $d\in X$, we denote the composite map
\[
f[gb,d]\xrightarrow H [fg b,fd]\xrightarrow u [b,fd]
\]
by $\vartheta$.

\begin{lemma}\label{delta.thm}
Let $((?)^*,(?)_*)$ be an adjoint pair where $(?)_*$ is a covariant
monoidal almost pseudofunctor on a category $\Cal S$ and $X_*$ is a 
symmetric monoidal closed category for $X\in \Cal S$.
Then for morphisms $f:X\rightarrow Y$ and $g:Y\rightarrow Z$ of $\Cal S$
and $b,d\in Z^*$, the diagram
\begin{equation}\label{delta.eq}
\xymatrix{
(gf)^*(b\otimes d) \ar[rr]^\Delta \ar[d]^{d^{-1}} & &
(gf)^*b\otimes (gf)^*d \ar[d]^{d^{-1}\otimes d^{-1}} \\
f^*g^*(b\otimes d) \ar[r]^-\Delta &
f^*(g^*b\otimes g^*d) \ar[r]^\Delta &
f^*g^*b\otimes f^*g^*d
}
\end{equation}
is commutative.
\end{lemma}

\proof Consider the diagram
\begin{landscape}
\[
\xymatrix{
(gf)^*(b\otimes d) \ar[r]^-{u\otimes u} \ar[d]^-{d^{-1}} &
(gf)^*((gf)_*(gf)^*b\otimes (gf)_*(gf)^*d)
\ar[r]^-m \ar[d]^{d^{-1}} &
(gf)^*(gf)_*((gf)^*b\otimes (gf)^*d)
\ar[r]^-\varepsilon \ar[ddd]^{d^{-1}c} &
(gf)^*b\otimes(gf)^*d \ar[dddd]^{d^{-1}} \\
f^*g^*(b\otimes d) \ar[r]^-{u\otimes u} \ar[d]^{u\otimes u}
\ar@{}[dr]|{\text{(a)}} &
f^*g^*((gf)_*(gf)^*b\otimes (gf)_*(gf)^*d)
\ar[d]^{d^{-1}c}
\ar@{}[dr]|{\text{(b)}} & 
\ar@{}[dr]|{\text{(c)}} & \\
f^*g^*(g_*g^*b\otimes g_*g^*d)
\ar[r]^-{u\otimes u} \ar[d]^m &
f^*g^*(g_*f_*f^*g^*b\otimes g_*f_*f^*g^*d)
\ar[d]^m & & \\
f^*g^*g_*(g^*b\otimes g^*d)\ar[r]^-{u\otimes u} \ar[d]^\varepsilon &
f^*g^*g_*(f_*f^*g^*b\otimes f_*f^*g^*d)\ar[r]^-{m} \ar[d]^\varepsilon &
f^*g^*g_*f_*(f^*g^*b\otimes f^*g^*d) \ar[d]^\varepsilon & \\
f^*(g^*b\otimes g^*d) \ar[r]^-{u\otimes u} &
f^*(f_*f^*g^*b\otimes f_*f^*g^*d) \ar[r]^-{m} &
f^*f_*(f^*g^*b\otimes f^*g^*d) \ar[r]^-{\varepsilon} &
f^*g^*b\otimes f^*g^*d
}.
\]
\end{landscape}
(a) is commutative by \cite[Lemma~1.13]{ETI}.
The commutativity of (b) is one of our assumptions, see 
\cite[(3.6.7.2)]{Lipman}.
(c) is commutative by \cite[Lemma~1.14]{ETI}.
Commutativity of the other squares is trivial.
Thus the whole diagram is commutative, and we are done.
\qed

\begin{lemma}\label{Delta-commutativity.thm}
Let $f:X\rightarrow Y$ be a symmetric monoidal functor between symmetric
monoidal categories, and $g:Y\rightarrow X$ its adjoint.
For $b\in X$ and $d\in Y$, the diagram
\[
\xymatrix{
g(fb\otimes d) \ar[rr]^\Delta \ar[d]^u & &
gfb\otimes gd \ar[d]^\varepsilon \\
g(fb\otimes fgd) \ar[r]^m &
gf(b\otimes gd) \ar[r]^\varepsilon &
b\otimes gd
}
\]
is commutative.
\end{lemma}

\proof Follows from the commutativity of the diagram
\[
\xymatrix{
g(fb\otimes d) \ar[r]^-{u\otimes u} \ar[d]^u & 
g(fgfb \otimes fgd) \ar[r]^m \ar[d]^\varepsilon &
gf(gfb\otimes gd) \ar[r]^-\varepsilon \ar[d]^\varepsilon &
gfb\otimes gd \ar[d]^\varepsilon \\
g(fb\otimes fgd) \ar[ur]^u \ar[r]^\id &
g(fb\otimes fgd) \ar[r]^m &
gf(b\otimes gd) \ar[r]^\varepsilon &
b\otimes gd
}.\qed
\]

\begin{lemma}\label{vartheta-isom.thm}
Viewed as a functor on $?$, $\vartheta:f[gb,?]\rightarrow [b,?]f$ is 
right conjugate to $\Delta:g(?\otimes b)\rightarrow g?\otimes gb$.
In particular, if $(f,g)$ is a Lipman symmetric monoidal adjoint pair,
then $\vartheta$ is an isomorphism.
\end{lemma}

\proof Follows from the commutativity of the diagram
\begin{landscape}
\[
\xymatrix{
f[gb,d] \ar[r]^-{\trace} \ar[d]^{\trace} \ar@{}[dr]|{\text{(a)}} &
[b,f[gb,d]\otimes b] \ar[r]^-{u} \ar[d]^{u} &
[b,fg(f[gb,d]\otimes b)] \ar[r]^-{\Delta} \ar[d]^{u} &
[b,f(gf[gb,d]\otimes gb)] \ar[dd]^\varepsilon \\
[fgb,f[gb,d]\otimes fgb] \ar[r]^-{u} \ar[d]^m &
[b,f[gb,d]\otimes fgb] \ar[r]^-{u} \ar[d]^m &
[b,fg(f[gb,d]\otimes fg b)] \ar[d]^m
\ar@{}[r]|{\text{~~~(b)}} & \\
[fgb,f([gb,d]\otimes gb)] \ar[r]^-u \ar[d]^\ev &
[b,f([gb,d]\otimes gb)] \ar[r]^-u \ar[d]^\ev &
[b,fgf([gb,d]\otimes gb)] \ar[r]^-\varepsilon \ar[d]^\ev &
[b,f([gb,d]\otimes gb)] \ar[d]^{\ev} \\
[fgb,fd] \ar[r]^-{u} &
[b,fd] \ar[r]^-{u} \ar `d[drr] `r[rr]^\id [rr] &
[b,fgfd] \ar[r]^-{\varepsilon} &
[b,fd]\\
 & & & 
},
\]
\end{landscape}
where the commutativity of (a) and (b) follows from \cite[(1.32)]{ETI} and
Lemma~\ref{Delta-commutativity.thm}, respectively.
\qed

Consider that the diagram (\ref{delta.eq}) is that of functors on $b$
(consider that $d$ is fixed), and then take a conjugate diagram, we
immediately have:

\begin{lemma}
Let $\Cal S$, $((?)^*,(?)_*)$, $f$, and $g$ be as in {\rm
Lemma~\ref{delta.thm}}.
Then for $d\in Z_*$ and $e\in X_*$, the diagram
\[
\xymatrix{
[d,(gf)_*e] & & (gf)_*[(gf)^*d,e] \ar[ll]_{\vartheta} \\
[d,g_*f_*e] \ar[u]_{c^{-1}} &
g_*[g^*d,f_*e] \ar[l]_\vartheta &
g_*f_*[f^*g^*d,e] \ar[l]_\vartheta \ar[u]_{c^{-1}d^{-1}}
}
\]
is commutative.
\qed
\end{lemma}

\begin{lemma}\label{H-vartheta.thm}
Let $\Cal S$ and $((?)^*,(?)_*)$ be as in {\rm
Lemma~\ref{delta.thm}}.
Let 
\[
\xymatrix{
X'\ar[r]^{f'} \ar[d]^{g'} & Y' \ar[d]^g \\
X \ar[r]^f & Y
}
\] 
be a commutative diagram in $\Cal S$.
Then for $b\in X_*$ and $d\in X'_*$, the diagram
\[
\xymatrix{
f_*g'_*[(g')^*b,d] \ar[r]^{\vartheta} \ar[d]^c &
f_*[b,g'_*d] \ar[r]^H &
[f_*b,f_*g'_*d] \ar[dr]^c & \\
g_*f'_*[(g')^*b,d] \ar[r]^H &
g_*[f'_*(g')^*b,f'_*d] \ar[r]^\theta &
g_*[g^*f_*b,f'_*d] \ar[r]^\vartheta &
[f_*b,g_*f'_*d]
}
\]
is commutative.
\end{lemma}

\proof Consider the diagram
\[
\xymatrix{
f_*g'_*[(g')^*b,d] \ar[r]^-H \ar[d]^c &
f_*[g'_*(g')^*b,g'_*d] \ar[r]^-u \ar[d]^H 
\ar@{}[dr]|{\text{(b)}} &
f_*[b,g'_*d] \ar[d]^H \\
g_*f'_*[(g')^*b,d] \ar@{}[r]|{\text{(a)~~}} \ar[d]^H &
[f_*g'_*(g')^*b,f_*g'_*d] \ar[d]^c \ar[r]^-u &
[f_*b,f_*g'_*d] \ar[dd]^c \\
g_*[f'_*(g')^*b,f'_*d] \ar@{}[dr]|{\text{(c)}} \ar[d]^\theta \ar[r]^-H &
[g_*f'_*(g')^*b,g_*f'_*d] \ar@{}[r]|{\text{~~~~(d)}} \ar[d]^\theta & \\
g_*[g^*f_*b,f'_*d] \ar[r]^-H &
[g_*g^*f_*b,g_*f'_*d] \ar[r]^-u &
[f_*b,g_*f'_*d]
}.
\]
(a) is commutative by \cite[Lemma~1.39]{ETI}.
The commutativity of (b) and (c) is trivial.
(d) is commutative by \cite[Lemma~1.24]{ETI}.
\qed

\paragraph Let $f:X\rightarrow Y$ be a Lipman 
symmetric monoidal functor between
closed categories, and $g:Y\rightarrow X$ its adjoint.
We denote the composite
\[
fg[b,d]
\xrightarrow{P}
f[gb,gd]
\xrightarrow{\vartheta}
[b,fgd]
\]
by $\fG$.

\begin{lemma}
Let $\Cal S$, $((?)^*,(?)_*)$ and
\[
\xymatrix{
X'\ar[r]^{f'} \ar[d]^{g'} & Y' \ar[d]^g \\
X \ar[r]^f & Y
}
\] 
be as in {\rm Lemma~\ref{H-vartheta.thm}}.
Then for $b,d\in X_*$, the diagram
\[
\xymatrix{
g_*g^*f_*[b,d] \ar[r]^H \ar[d]^{c\theta} &
g_*g^*[f_*b,f_*d] \ar[r]^{\fG} &
[f_*b,g_*g^*f_*d] \ar[d]^{c\theta} \\
f_*g'_*(g')^*[b,d] \ar[r]^{\fG} &
f_*[b,g'_*(g')^*d] \ar[r]^H &
[f_*b,f_*g'_*(g')^*d]
}
\]
is commutative.
\end{lemma}

\proof Left to the reader.
Use \cite[(1.24)]{ETI}, \cite[(1.39)]{ETI}, and \cite[(1.59)]{ETI}.
\qed

\begin{lemma}
Let $\Cal S$, $((?)^*,(?)_*)$ and
\[
\xymatrix{
X'\ar[r]^{f'} \ar[d]^{g'} & Y' \ar[d]^g \\
X \ar[r]^f & Y
}
\] 
be as in {\rm Lemma~\ref{H-vartheta.thm}}.
Assume that $((?)^*,(?)_*)$ is Lipman.
Then for $b,d\in Y_*$,
the diagram
\[
\xymatrix{
g^*f_*f^*[b,d] \ar[r]^-\fG \ar[d]^{d\theta} &
g^*[b,f_*f^*d] \ar[r]^-P &
[g^*b,g^*f_*f^*d] \ar[d]^{d\theta} \\
f'_*(f')^*g^*[b,d] \ar[r]^-P &
f'_*(f')^*[g^*b,g^*d] \ar[r]^-{\fG} &
[g^*b,f'_*(f')^*g^*d]
}
\]
is commutative.
\end{lemma}

\proof Left to the reader.
Use \cite[(1.26)]{ETI}, \cite[(1.54)]{ETI}, and \cite[(1.59)]{ETI}.
\qed

\begin{lemma}\label{trivial-isom.thm}
Let $I$ be a small category, and $f:X\rightarrow Y$ a morphism of 
$I\op$-diagrams of schemes.
Then for $\M\in\Mod(Y)$ and $\N\in\Mod(X)$, the composite
\[
\vartheta: f_*\uHom_{\O_X}(f^*\M,\N)
\xrightarrow{H}
\uHom_{\O_Y}(f_*f^*\M,f_*\N)
\xrightarrow{u}
R\uHom_{\O_Y}(\M,f_*\N)
\]
is an isomorphism.
\end{lemma}

\proof This is an immediate consequence of Lemma~\ref{vartheta-isom.thm}.
\qed

\begin{lemma}\label{cartesian-immersion.thm}
Let $I$ be a small category, and $f:U\rightarrow X$ be an
open immersion of $I\op$-diagrams of schemes.
Let $\M,\N\in\Mod(X)$.
If either
\begin{description}
\item[(i)] $\M$ is equivariant; or
\item[(ii)] $f$ is cartesian,
\end{description}
Then the canonical map
\[
P:f^*\uHom_{\O_X}(\M,\N)\rightarrow\uHom_{\O_U}(f^*\M,f^*\N)
\]
is an isomorphism of presheaves.
In particular, it is an isomorphism of sheaves.
\end{lemma}

\proof {\bf (ii)}
Taking the section at $(i,V)$, where $i\in I$ and $V\in\Zar(U_i)$, 
it suffices to show that the map induced by the restriction
\begin{equation}\label{restriction.eq}
\Hom_{\Zar(X)/(i,V)}(\M|_{(i,V)},\N|_{(i,V)})
\rightarrow
\Hom_{\Zar(U)/(i,V)}(\M|_{(i,V)},\N|_{(i,V)})
\end{equation}
is an isomorphism, see the description of $P$ in \cite[(2.8)]{HO}.
But as $U$ is cartesian, $\Zar(U)/(i,V)\hookrightarrow \Zar(X)/(i,V)$
is an equivalence.
Indeed, if $(j,W)\rightarrow (i,V)$ is a morphism in $\Zar(X)$, it must
be a morphism in $\Zar(U)$.
Thus (\ref{restriction.eq}) is an isomorphism, and we are done.

{\bf (i)} Similarly to the proof of \cite[(2.13)]{HO}, the problem is
reduced to the case of single schemes.
Then the assertion follows from {\bf (ii)} immediately.
\qed

\begin{lemma}\label{commutative-G.thm}
Let $((?)^*,(?)_*)$ be a Lipman 
monoidal adjoint pair on a category $\Cal S$ where $X_*$ is
closed for every $X\in \Cal S$.
For morphisms $g:X\rightarrow Y$ and $f:Y\rightarrow Z$ 
of $\Cal S$ and $a,b\in Z^*$, 
the composite
\[
f_*f^*[a,b]
\xrightarrow{\fG}
[a,f_*f^*b]
\xrightarrow u
[a,f_*g_*g^*f^*b]
\]
agrees with the composite
\begin{multline*}
f_*f^*[a,b]
\xrightarrow u
f_*g_*g^*f^*[a,b]
\xrightarrow{dc^{-1}}
(fg)_*(fg)^*[a,b]
\xrightarrow{\fG}\\
[a,(fg)_*(fg)^*b]
\xrightarrow{d^{-1}c}
[a,f_*g_*g^*f^*b].
\end{multline*}
\end{lemma}

\proof Left to the reader.
Use \cite[(1.39), (1.54), (1.56)]{ETI}.
\qed

\begin{corollary}
Let $((?)^*,(?)_*)$ and $g:X\rightarrow Y$ be as in 
{\rm Lemma~\ref{commutative-G.thm}}.
Then the composite
\[
[a,b]\xrightarrow u
g_*g^*[a,b]
\xrightarrow{\fG}
[a,g_*g^*b]
\]
is $u$.
\end{corollary}

\proof Let $f=\id$ in Lemma~\ref{commutative-G.thm}.
\qed

\begin{lemma}\label{fG-isom.thm}
Let $I$ be a small category, $X$ an $I\op$-diagram of schemes,
$f:U\hookrightarrow X$ an open subdiagram.
Let $\M,\N\in \Mod(X)$, and consider the map
\[
\fG:f_*f^*\uHom_{\O_X}(\M,\N)\rightarrow
\uHom_{\O_X}(\M,f_*f^*\N).
\]
If $f$ is cartesian or $\M$ is equivariant, then $\fG$ is an isomorphism.
\end{lemma}

\proof Note that $\fG$ is the composite
\[
f_*f^*\uHom_{\O_X}(\M,\N)\xrightarrow{P}
f_*\uHom_{\O_U}(f^*\M,f^*\N)\xrightarrow{\vartheta}
\uHom_{\O_X}(\M,f_*f^*\N).
\]
$P$ is an isomorphism by Lemma~\ref{cartesian-immersion.thm}.
$\vartheta$ is an isomorphism by Lemma~\ref{trivial-isom.thm}.
So $\fG$ is an isomorphism.
\qed

\paragraph \label{fH-settings.par}
Let $I$ be a small category, $X$ an $I\op$-diagram of schemes,
$f:U\hookrightarrow X$ an open subdiagram, and $g:V\hookrightarrow U$ an
open subdiagram.
Then for $\M,\N\in\Mod(X)$, 
the diagram
\begin{equation}\label{fH-def.eq}
\xymatrix{
0 \ar[d] & 0 \ar[d] \\
\uGamma_{U,V}\uHom_{\O_X}(\M,\N) \ar[d]^-\iota &
\uHom_{\O_X}(\M,\uGamma_{U,V}\N) \ar[d]^-\iota \\
f_*f^*\uHom_{\O_X}(\M,\N) \ar[d]^-u \ar[r]^-{\fG} &
\uHom_{\O_X}(\M,f_*f^*\N)\ar[d]^-u \\
f_*g_*g^*f^*\uHom_{\O_X}(\M,\N)~~~ \ar[r]^{d^{-1}c\fG dc^{-1}} &
~~~\uHom_{\O_X}(\M,f_*g_*g^*f^*\N)
}
\end{equation}
is commutative with exact columns by
Lemma~\ref{commutative-G.thm}.
So there is a unique natural map
\begin{equation}\label{fH.eq}
\fH: \uGamma_{U,V}\uHom_{\O_X}(\M,\N)\rightarrow 
\uHom_{\O_X}(\M,\uGamma_{U,V}\N)
\end{equation}
such that $\iota\fH=\fG\iota$.

\begin{lemma}\label{fH-isom.thm}
Let the notation be as in {\rm(\ref{fH-settings.par})}.
If both $f$ and $g$ are cartesian, or $\M$ is equivariant, 
then $\fH$ in {\rm(\ref{fH.eq})} is an isomorphism.
\end{lemma}

\proof Follows from Lemma~\ref{fG-isom.thm} and the five lemma applied to
the diagram (\ref{fH-def.eq}).
\qed

\paragraph \label{preserving-injectives.par}
Let $I$ be a small category, and $f:U\rightarrow X$ be an
open immersion of $I\op$-diagram of schemes.
Then $\Gamma((i,V),f^*_\hst\M)=\Gamma((i,V),\M)$ for $\M\in\hst(X)$ 
almost by definition,
where $\hst=\PM$ or $\Mod$.
Thus if $j:\Zar U\hookrightarrow  \Zar X$ is the inclusion,
then $f^*_{\hst}=j^\#_{\hst}$.
Thus $f^*_\hst$ has a left adjoint $j_\#^\hst$, 
as well as the right adjoint $f_*$.
Hence $f^*_\hst$ preserves arbitrary limits as well as arbitrary colimits.
We denote $j_\#^{\hst}$ by $f_!$ or $f_!^{\hst}$ by an obvious reason.

Note that $\Gamma((i,V),f_!^{\PM}(\M))$ is 
$\Gamma((i,V),\M)$ if $V\subset U_i$,
and zero if $V\not\subset U_i$.
In particular, $f_!^{\PM}$ is exact.

Note also that we have a commutative diagram
\[
\xymatrix{
\Zar(U_i) \ar[r]^{Q(U,i)} \ar[d]^{j} & \Zar(U) \ar[d]^j \\
\Zar(X_i) \ar[r]^{Q(X,i)} & \Zar(X)
},
\]
where $Q(X,i)$ and $Q(U,i)$ are obvious inclusions, see \cite[(4.5)]{ETI}.
By \cite[(2.57)]{ETI}, Lipman's theta \cite[(1.21)]{ETI} 
$\theta:j_\#^{\PM} Q(U,i)^\#\rightarrow Q(X,i)^\# j_\#^{\PM}$, namely, 
$\theta: (f_i)_!^{\PM}(?)_i\rightarrow (?)_i f_!^{\PM}$ at $(i,V)$ is
the identity of $\Gamma((i,V),\M)$ if $V\subset U_i$, and zero otherwise.
In particular, $\theta$ is an isomorphism.

Note that $f_!^{\Mod}=j_\#^{\Mod}=aj_{\#}^{\PM}q=af_!^{\PM}q$.
By \cite[(2.59)]{ETI}, 
$\theta:(?)_if_!^{\Mod}\rightarrow (f_i)_!^{\Mod}(?)_i$ 
is an isomorphism.
It is well-known that $(f_i)_!^{\Mod}
$ is exact, and hence $f_!^{\Mod}$ is exact.

Since $f^*_{\hst}$ has an exact left adjoint $f_!^{\hst}$, 
$f^*_{\hst}$ preserves injectives and $K$-injectives for $\hst=\PM,\Mod$.

\begin{lemma}\label{limits.thm}
Let the notation be as in {\rm (\ref{fH-settings.par})}.
Then $f^*$, $(fg)^*$, and $\uGamma_{U,V}$ preserves arbitrary limits.
\end{lemma}

\proof 
By the discussion in (\ref{preserving-injectives.par}), 
$f^*$ and $(fg)^*$ preserves limits.

Now let $(\M_\lambda)$ be a system in $\Mod(X)$.
Then
\[
\xymatrix{
0 \ar[r] &
\uGamma_{U,V}\projlim\M_\lambda \ar[d] \ar[r]^\iota &
f_*f^*\projlim\M_\lambda \ar[d]^{\cong} \ar[r]^-u &
f_*g_*g^*f^*\projlim \M_\lambda \ar[d]^{\cong}\\
0 \ar[r] &
\projlim \uGamma_{U,V}\M_\lambda \ar[r]^\iota &
\projlim f_*f^*\M_\lambda \ar[r]^-u &
\projlim f_*g_*g^*f^*\M_\lambda
}
\]
is a commutative diagram with exact rows.
By the five lemma, $\uGamma_{U,V}$ preserves limits.
\qed

\paragraph Let the notation be as in (\ref{fH-settings.par}).
For a complex $\Bbb F$ in $\Mod(X)$, a natural map
\[
\fH: \uGamma_{U,V}\uHom_{\O_X}(\Bbb F,?)\rightarrow
\uHom_{\O_X}(\Bbb F,?)\uGamma_{U,V}
\]
between functors on the category of complexes in $\Mod(X)$ is defined.
By Lemma~\ref{limits.thm} and Lemma~\ref{fH-isom.thm},
it is an isomorphism if $f$ and $g$ are cartesian, or $\Bbb F$ is a
complex of equivariant sheaves.
Similarly, 
\[
\fG: 
f_*f^*\uHom_{\O_X}(\Bbb F,?)\rightarrow
\uHom_{\O_X}(\Bbb F,?)f_*f^*
\]
and
\[
d^{-1}c\fG dc^{-1}:
f_*g_*g^*f^*\uHom_{\O_X}(\Bbb F,?)\rightarrow
\uHom_{\O_X}(\Bbb F,?)f_*g_*g^*f^*
\]
are induced.

\begin{lemma} Let $(\Bbb X,\O_{\Bbb X})$ be a ringed site,
$\Bbb F$ a complex of $\O_{\Bbb X}$-modules, and $\Bbb G$ a $K$-injective
complex of $\O_{\Bbb X}$-modules.
Then $\uHom_{\O_{\Bbb X}}(\Bbb F,\Bbb G)$ is weakly $K$-injective.
\end{lemma}

\proof Let $\Bbb H$ be any exact $K$-flat complex.
Then
\[
\Hom_{\O_{\Bbb X}}(\Bbb H,\uHom_{\O_{\Bbb X}}(\Bbb F,\Bbb G))
\cong
\Hom_{\O_{\Bbb X}}(\Bbb H\otimes \Bbb F,\Bbb G)
\]
is exact, since 
$\Bbb H\otimes \Bbb F$ is exact \cite[Lemma~3.21, {\bf 2}]{ETI}
and $\Bbb G$ is $K$-injective.
By \cite[Lemma~3.25, {\bf 5}]{ETI}, $\uHom_{O_{\Bbb X}}(\Bbb F,\Bbb G)$
is weakly $K$-injective.
\qed

\begin{lemma}
The canonical maps
\[
\zeta:R(\uGamma_{U,V}\uHom_{\O_X}(\Bbb F,?))
\rightarrow R\uGamma_{U,V} R\uHom_{\O_X}(\Bbb F,?),
\]
\[
\zeta:R(f_*f^*\uHom_{\O_X}(\Bbb F,?))
\rightarrow Rf_*f^* R\uHom_{\O_X}(\Bbb F,?),
\]
and
\[
\zeta:R(f_*g_*g^*f^*\uHom_{\O_X}(\Bbb F,?))
\rightarrow Rf_*Rg_*g^*f^* R\uHom_{\O_X}(\Bbb F,?)
\]
are isomorphisms.
\end{lemma}

\proof 
For a $K$-injective complex $\Bbb G$, $\uHom_{\O_X}(\Bbb F,\Bbb G)$ is
weakly $K$-injective.
So it is $K$-flabby, and $\uGamma_{U,V}$-acyclic \cite[(4.3)]{HO}.
In particular, $f^*\uHom_{\O_X}(\Bbb F,\Bbb G)$ and
$g^*f^*\uHom_{\O_X}(\Bbb F,\Bbb G)$ are $K$-limp by \cite[(4.6)]{HO},
and the assertion follows.
\qed

\paragraph By the lemma,
the composite
\begin{multline*}
\fH: R\uGamma_{U,V}R\uHom_{\O_X}(\Bbb F,?)
\xrightarrow{\zeta^{-1}}
R(\uGamma_{U,V}\uHom_{\O_X}(\Bbb F,?))\\
\xrightarrow{R\fH}
R(\uHom_{\O_X}(\Bbb F,?)\uGamma_{U,V})
\xrightarrow{\zeta}
R\uHom_{\O_X}(\Bbb F,?)R\uGamma_{U,V}
\end{multline*}
is defined.
Similarly,
\[
\fG: 
Rf_*f^*R\uHom_{\O_X}(\Bbb F,?)
\rightarrow
R\uHom_{\O_X}(\Bbb F,?)Rf_*f^*
\]
and
\[
d^{-1}c\fG dc^{-1}:
Rf_*Rg_*g^*f^*R\uHom_{\O_X}(\Bbb F,?)
\rightarrow
R\uHom_{\O_X}(\Bbb F,?)Rf_*Rg_*g^*f^*
\]
are induced.
Note that
\begin{equation}\label{triangle.eq}
\xymatrix{
R\uGamma_{U,V}R\uHom_{\O_X}(\Bbb F,?)
\ar[r]^{\fH}
\ar[d]^\iota &
R\uHom_{\O_X}(\Bbb F,?)R\uGamma_{U,V}
\ar[d]^\iota \\
Rf_*f^*R\uHom_{\O_X}(\Bbb F,?)
\ar[r]^{\fG}
\ar[d]^u &
R\uHom_{\O_X}(\Bbb F,?)Rf_*f^*
\ar[d]^u \\
Rf_*Rg_*g^*f^*R\uHom_{\O_X}(\Bbb F,?)~~~
\ar[r]^{d^{-1}c\fG dc^{-1}} \ar[d] &
~~~R\uHom_{\O_X}(\Bbb F,?)Rf_*Rg_*g^*f^* \ar[d] \\
R\uGamma_{U,V}R\uHom_{\O_X}(\Bbb F,?)[1] \ar[r]^{\fH[1]} &
R\uHom_{\O_X}(\Bbb F,?)R\uGamma_{U,V}[1]
}
\end{equation}
is a commutative diagram with columns being triangles.

\begin{lemma}\label{derived-trivial.thm}
Let $I$ be a small category, and $f:X\rightarrow Y$ a morphism of 
$I\op$-diagrams of schemes.
Then the composite
\begin{multline*}
\vartheta: Rf_*R\uHom_{\O_X}(Lf^*\Bbb F,\Bbb G)
\xrightarrow{H}
R\uHom_{\O_Y}(Rf_*Lf^*\Bbb F,Rf_*\Bbb G)\\
\xrightarrow{u}
R\uHom_{\O_Y}(\Bbb F,Rf_*\Bbb G)
\end{multline*}
is an isomorphism between functors on $D(Y)\op\times D(X)$.
\end{lemma}

\proof This is an immediate consequence of \cite[(1.49)]{ETI} and
\cite[(8.23), {\bf 5}]{ETI}.
\qed

\begin{corollary}\label{cartesian-derived.thm}
Let $f:U\rightarrow X$ be a cartesian open immersion.
Then
\[
P:f^*R\uHom_{\O_X}(\Bbb F,\Bbb G)\rightarrow R\uHom_{\O_U}(f^*\Bbb F,f^*
\Bbb G)
\]
is an isomorphism for any $\Bbb F,\Bbb G\in D(X)$.
\end{corollary}

\proof If $\Bbb G$ is a $K$-injective complex in $K(X)$, then 
so is $f^*\Bbb G$ by (\ref{preserving-injectives.par}).
So it suffices to show that
\[
f^*\uHom_{\O_X}(\Bbb F,\Bbb G)\rightarrow \uHom_{\O_U}(f^*\Bbb F,f^*\Bbb G)
\]
is an isomorphism of complexes, 
if $\Bbb F$ and $\Bbb G$ are complexes in $\Mod(X)$.
This follows from Lemma~\ref{cartesian-immersion.thm} and the fact
that $f^*$ preserves direct product.
\qed

\begin{lemma}\label{derived-P.thm}
Let $I$ be a small category, and $f:X\rightarrow Y$ a morphism of 
$I\op$-diagrams of schemes.
Let $\Bbb F$ and $\Bbb G$ be objects in $D(Y)$.
Assume that one of the following holds:
\begin{description}
\item[(i)] $f$ is locally an open immersion, $\Bbb F\in D_{\EM}(Y)$,
and one of the following holds:
\begin{description}
\item[(a)] $\Bbb G\in D^+(Y)${\rm;}
\item[(b)] $\Bbb F\in D^+_{\EM}(Y)${\rm;}
\item[(c)] $\Bbb G\in D_{\Lqc}(Y)$.
\end{description}
\item[(ii)] $f$ is flat, $Y$ is locally noetherian, $\Bbb G\in D^+(Y)$,
and $\Bbb F\in D^-_{\Coh}(Y)$.
\item[(iii)] $f$ is flat, $Y$ is locally noetherian, 
$\Bbb F\in D_{\Coh}(Y)$, and both $\Bbb G$ and $f^*\Bbb G$ have finite
injective dimension.
\end{description}
Then the canonical map
\[
P: f^*R\uHom_{\O_Y}(\Bbb F,\Bbb G)
\rightarrow
R\uHom_{\O_X}(f^*\Bbb F,f^*\Bbb G)
\]
is an isomorphism.
\end{lemma}

\proof Similarly to \cite[Lemma~1.59]{ETI}, using \cite[Lemma~1.56]{ETI},
it is easy to prove that 
the diagram
\[
\xymatrix{
(?)_if^*R\uHom_{\O_Y}(\Bbb F,\Bbb G)\ar[r]^P \ar[d]^{\theta^{-1}} &
(?)_iR\uHom_{\O_X}(f^*\Bbb F,f^*\Bbb G) \ar[d]^H \\
f_i^*(?)_i R\uHom_{\O_Y}(\Bbb F,\Bbb G) \ar[d]^H&
R\uHom_{\O_{X_i}}((?)_if^*\Bbb F,(?)_if^*\Bbb G) 
\ar[d]^{[\theta,\theta^{-1}]}\\
f_i^*R\uHom_{\O_{Y_i}}(\Bbb F_i,\Bbb G_i)
\ar[r]^P &
R\uHom_{\O_{X_i}}(f_i^*\Bbb F_i,f_i^*\Bbb G_i)
}
\]
is commutative for $i\in I$.
Note that the vertical morphisms are isomorphisms by \cite[(13.9)]{ETI} and
\cite[(6.25)]{ETI}.
So in order to prove that the top $P$ is an isomorphism for each $i\in I$,
it suffices to prove the bottom $P$ is an isomorphism.
So we may assume that the problem is on single schemes.

{\bf (i)} We may assume that $f$ is an open immersion.
Then this is a special case of Corollary~\ref{cartesian-derived.thm}.

{\bf (ii)} This is \cite[(5.8)]{Hartshorne}.

{\bf (iii)} This follows from {\bf (ii)} and the way-out lemma
(\cite[Proposition~I.7.1, (iii)]{Hartshorne}).
\qed

\begin{theorem}\label{fH-isom-derived.thm}
Let $I$ be a small category, $X$ an $I\op$-diagram of schemes,
$f:U\hookrightarrow X$ an open subdiagram,
and $g:V\hookrightarrow U$ an open subdiagram.
Let $\Bbb F$ and $\Bbb G$ be in $D(X)$.
If one of the following holds, then 
\[
\fH:R\uGamma_{U,V}R\uHom_{\O_X}(\Bbb F,?)
\rightarrow
R\uHom_{\O_X}(\Bbb F,?)R\uGamma_{U,V}
\]
is an isomorphism:
\begin{description}
\item[(i)] $f$ and $g$ are cartesian;
\item[(ii)] $\Bbb F\in D_{\EM}(X)$, and one of the following hold:
{\rm(a)} $\Bbb G\in D^+(X)$; {\rm(b)} $\Bbb F\in D^+_{\EM}(X)$;
{\rm(c)} $\Bbb G\in D_{\Lqc}(X)$.
\end{description}
\end{theorem}

\proof By Lemma~\ref{derived-trivial.thm} and
Lemma~\ref{derived-P.thm}, the two maps 
$\vartheta P$ and $d^{-1}c\vartheta Pdc^{-1}$ in 
(\ref{triangle.eq}) are isomorphisms.
As the columns of (\ref{triangle.eq}) are triangles,
the third horizontal map $\fH$ is also an isomorphism.
\qed

\section{Matlis duality and the local duality}

Let $S$ be a scheme, $G$ an $S$-group scheme, $(X,Y)$ a standard
$G$-local $G$-scheme.
That is, $X$ is a standard
$G$-local $G$-scheme, and $Y$ is its unique minimal
closed $G$-subscheme.
We denote the inclusion $Y\hookrightarrow X$ by $j$.

We denote the defining ideal sheaf of $Y$ by $\I$.
Thus $\I$ is the unique $G$-maximal $G$-ideal of $\O_X$.
We fix the generic point of an irreducible component of $Y$ and
denote it by $\eta$.

\begin{lemma}\label{m-r-local.thm}
Let $\C$ be a class of noetherian local rings.
Assume that if $A\in\C$ and $B$ is essentially of finite type over $A$,
then $B\in\C$.
Let $\Bbb P(A,M)$ be a property of a pair $(A,M)$ of a 
finitely generated module $M$ over a noetherian local ring $A$ such that
$A\in\C$.
Assume that
\begin{description}
\item[(i)] If $A\in\C$,
$\Bbb P(A,M)$ holds, and $P\in \Spec A$, then $\Bbb P(A_P,M_P)$
holds.
\item[(ii)] If $A\in\C$, $M$ a finite $A$-module,
and $A\rightarrow B$ is a flat local
homomorphism essentially of finite type with local complete intersection
fibers \(resp.\ geometrically regular fibers\), then $\Bbb P(A,M)$ holds
if and only if $\Bbb P(B,B\otimes_AM)$ holds.
\end{description}
Assume that the all local rings of $X$ belong to $\C$.
For $\M\in\Coh(G,X)$, if $\Bbb P(\O_{X,\eta},\M_\eta)$ holds
\(resp.\ $\Bbb P(\O_{X,\eta},\M_\eta)$ holds and 
either the second projection $p_2:G\times X\rightarrow X$ is smooth or
$S=\Spec k $ with $k$ a perfect field and $G$ is of finite type over $S$\), 
then
$\Bbb P(\O_{X,x},\M_x)$ holds for any $x\in X$.
\end{lemma}

\proof Let $Z$ be the unique integral
closed subscheme of $X$ whose generic point is
$x$.
Let $Z^*$ be the unique minimal closed $G$-subscheme of $X$ containing $Z$, 
see \cite{HM}.
As $\eta\in Y\subset Z^*$, there exists some irreducible component $Z_0
$ of $Z$ such that $\eta\in Z_0$.
Let $\zeta$ be the generic point of $Z_0$.
Since $\Bbb P(\O_{X,\eta},\M_\eta)$ holds and
$\zeta$ is a generalization of $\eta$, 
$\Bbb P(\O_{X,\zeta},\M_\zeta)$ holds.
Then by \cite[Corollary~7.6]{HM}, $\Bbb P(\O_{X,x},\M_x)$ holds.
\qed

\begin{corollary}\label{G-local-M-R.thm}
Let $m$, $n$, and $g$ be non-negative integers or $\infty$.
Then 
\begin{description}
\item[(i)]
Let $\M\in\Coh(G,X)$, and assume that $\M_\eta$ is maximal Cohen--Macaulay
\(resp.\ of finite injective dimension, projective dimension $m$, 
$\dim-\depth=n$, torsionless, reflexive, G-dimension $g$, zero\) as an
$\O_{X,\eta}$-module.
Then $\M_x$ is so as an $\O_{X,x}$-module for any $x\in X$.
\item[(ii)] If $\O_{X,\eta}$ is a complete intersection, 
then $X$ is locally a complete intersection.
\item[(iii)] Assume that $p_2:G\times X\rightarrow X$ is smooth or
$S=\Spec k$ with $k$ a perfect field and $G$ is of finite type over
$S$.
If $\O_{X,\eta}$ is regular, then $X$ is regular.
\item[(iv)] In addition to the assumption of {\bf (iii)},
assume further that $X$ is a locally excellent $\Bbb F_p$-scheme,
where $p$ is a prime number.
If $\O_{X,\eta}$ is $F$-regular \(resp.\ $F$-rational\), then the all
local rings of $X$ is $F$-regular \(resp.\ $F$-rational\).
\end{description}
\end{corollary}

\proof {\bf (i)} Let $\C$ be the class of all noetherian local rings,
and $\Bbb P(A,M)$ be ``$M$ is a maximal Cohen--Macaulay $A$-module.''
We can apply Lemma~\ref{m-r-local.thm}.
Similarly for other properties.

{\bf (ii)} Let $\C$ be the class of all noetherian local rings,
and $\Bbb P(A,M)$ be ``$A$ is a complete intersection.''
Then as $\Bbb P(\O_{X,\eta},0)$ holds, $\Bbb P(\O_{X,x},0)$ holds for
any $x\in X$.

{\bf (iii)} Let $\C$ be the class of all noetherian local rings,
and $\Bbb P(A,M)$ be ``$A$ is regular.''

{\bf (iv)} Let $\C$ be the class of all excellent noetherian local
rings of characteristic $p$, and $\Bbb P(A,M)$ be ``$A$ is $F$-regular''
(resp.\ ``$A$ is $F$-rational'').
\qed

\begin{corollary}\label{faithful-exact.thm}
The stalk functor $(?)_\eta:\Qch(G,X)\rightarrow \Mod(\O_{X,\eta})$ is
faithfully exact.
\end{corollary}

\proof The exactness is well-known.
Let $\M\in\Qch(G,X)$ and assume that $\M\neq 0$.
Then as $\Qch(G,X)$ is locally noetherian and its noetherian object is
nothing but a coherent $(G,\O_X)$-module, 
$\M$ contains a nonzero coherent $(G,\O_X)$-submodule $\N$.
Then by Corollary~\ref{G-local-M-R.thm}, $\M_\eta\supset \N_\eta\neq 0$.
This shows that $(?)_\eta$ is faithfully exact.
\qed

\begin{remark}
Formally, $(?)_\eta$ is a functor from $\Qch(G,X)=\Qch(B_G^M(X))$, or 
more generally, from $\Mod(G,X)=\Mod(B_G^M(X))$ to $\Mod(\O_{X,\eta})$, and
is the composite
\[
\Mod(B_G^M(X))\xrightarrow{(?)_0} \Mod(B_G^M(X)_0)=\Mod(X_0)
\xrightarrow{h^*}\Mod(\Spec \O_{X,\eta}),
\]
where $h:\Spec \O_{X,\eta}\rightarrow X_0$ is the inclusion.
Thus $(?)_\eta$ is sometimes written as $(?)_\eta(?)_0$,
where $(?)_\eta$ means $h^*$.
\end{remark}

\begin{corollary}[$G$-NAK]\label{G-NAK}
Let $\M\in\Coh(G,X)$.
If $j^*\M=0$, then $\M=0$, where $j:Y\hookrightarrow X$ is the inclusion.
\end{corollary}

\proof Since $j^*\M=0$, $\M/\I\M=0$.
So $\M_\eta/\I_\eta\M_\eta=0$.
By Nakayama's lemma, $\M_\eta=0$.
By Corollary~\ref{faithful-exact.thm}, $\M=0$.
\qed

\begin{proposition}\label{G-artinian-CM.thm}
A standard $G$-artinian $G$-scheme is Cohen--Macaulay.
\end{proposition}

\proof By Lemma~\ref{G-artinian.thm}, we may assume that the 
$G$-scheme is
$G$-local.
So let $X$ be a $G$-artinian $G$-local standard $G$-scheme.
Let $Y$, $\eta$, and $\I$ be as above.

Then $\grad{0}=\I$, since $\I$ is the only $G$-prime ideal
(for the definition and basic properties of $\grad{?}$, see 
\cite[section~4]{HM}).
So $Y=X$, set theoretically.
Thus $\eta$ is the generic point of an irreducible component of $X$.
So $\O_{X,\eta}$ is an artinian ring, and hence is Cohen--Macaulay.
By Corollary~\ref{G-local-M-R.thm}, $X$ is Cohen--Macaulay.
\qed

\begin{corollary} $Y$ is Cohen--Macaulay.
\end{corollary}

\proof Since $Y$ is $G$-artinian $G$-local standard, the corollary
follows immediately from 
Proposition~\ref{G-artinian-CM.thm}.
\qed

\paragraph From now on, we assume that $X$ has a $G$-dualizing complex 
$\Bbb I_X$ (see \cite[(31.15)]{ETI}).
For a $G$-morphism $f:X'\rightarrow X$ which is separated of finite type,
we denote $f^!\Bbb I_X$ by $\Bbb I_{X'}$, where $f^!$ is the 
twisted inverse functor $B^G_M(f)^!$ (see \cite[chapter~29]{ETI}).
Note that $\Bbb I_{X'}$ is a $G$-dualizing complex of $X'$ 
\cite[Lemma~31.11]{ETI}.
By \cite[Lemma~31.6]{ETI}, $\Bbb I_{X'}$, 
viewed as a complex of $\O_{X'}$-modules,
is a dualizing complex of $X'$.

Since $\O_{Y,\eta}$ is Cohen--Macaulay, 
there is only one $i$ such that
$H^i(\Bbb I_Y)_\eta\neq 0$.
This is equivalent to say that $H^i(\Bbb I_Y)\neq 0$.
If this $i$ is $0$, then we say that $\Bbb I_X$ is $G$-normalized.
If $X$ has a $G$-dualizing complex, then by shifting, $X$ has a 
$G$-normalized $G$-dualizing complex.

From now on, we always assume that $\Bbb I_X$ is $G$-normalized.

\begin{lemma}
$\Bbb I_{X,\eta}$ is a normalized dualizing complex of the local ring 
$\O_{X,\eta}$.
In particular, $H^0_{\fm_\eta}(\Bbb I_{X,\eta})$ is the injective hull of
the residue field $\kappa(\eta)$ of $\O_{X,\eta}$,
where $\fm_\eta$ is the maximal ideal of $\O_{X,\eta}$.
\end{lemma}

\proof Since $\Bbb I_X$ is a dualizing complex, $\Bbb I_{X,\eta}$ is
also a dualizing complex of $\O_{X,\eta}$.
We prove that $\Bbb I_{X,\eta}$ is normalized.
Let $\Bbb D$ be a normalized dualizing complex of $\O_{X,\eta}$,
and set $\Bbb I_{X,\eta}\cong \Bbb D[r]$.
We want to prove that $r=0$.

Consider the commutative diagram
\[
\xymatrix{
X & \Spec \O_{X,\eta} \ar[l]_-{p} \\
Y \ar[u]_{j} & \Spec \O_{Y,\eta} \ar[u]_{j'} \ar[l]_-{q}
}.
\]
By the commutativity with restrictions \cite[Proposition~18.14]{ETI},
\[
H^0(\Bbb I_Y)_\eta\cong 
H^0(q^*j^!\Bbb I_X)\cong H^0((j')^!\Bbb I_{X,\eta})
\cong \Ext^r_{\O_{X,\eta}}(\O_{Y,\eta},\Bbb D)\neq 0.
\]
The Matlis dual of the last module is $H^{-r}_{\fm_\eta}(\O_{Y,\eta})$, by
the local duality \cite[(V.6.3)]{Hartshorne}.
Since $\O_{Y,\eta}$ is an $\O_{X,\eta}$-module of finite length, 
$H^{-r}_{\fm_\eta}(\O_{Y,\eta})\neq 0$ implies $r=0$.
\qed

\begin{lemma}\label{E-def.thm}
$\uH^i_Y(\Bbb I_X)=0$ for $i\neq 0$, and 
$\uH^0_Y(\Bbb I_X)_\eta$ is the injective hull of the residue field 
$\kappa(\eta)$ of the local ring $\O_{X,\eta}$.
\end{lemma}

\proof By \cite[Theorem~6.10]{HO},
\[
(\uH^i_Y(\Bbb I_X))_\eta\cong H^i((?)_\eta R\uGamma_Y \Bbb I_X)
\cong H^i(R\uGamma_{\I_\eta}(?)_\eta\Bbb I_X)
\cong H^i_{\fm_\eta}(\Bbb I_{X,\eta}).
\]
Since $\Bbb I_{X,\eta}$ is a normalized dualizing complex of $\O_{X,\eta}$,
the last module is zero if $i\neq 0$ and is the injective hull of the
residue field $\kappa(\eta)$ of the local ring $\O_{X,\eta}$ if $i=0$.
As $(?)_\eta$ is faithfully exact, we are done.
\qed

\paragraph We set $\E:=\uH^0_Y(\Bbb I_X)$, and call it the 
{\em $G$-sheaf of Matlis}.
Note that the definition of $\E$ depends on the choice of $\Bbb I_X$.
Note also that $\E_\eta$ is the injective hull of the residue field of
$\O_{X,\eta}$.

\begin{lemma}\label{fid.thm}
$\E$ is of finite injective dimension as an object of 
$\Mod(G,X)$.
\end{lemma}

\proof We may assume that $\Bbb I_X$ is a bounded complex of injective objects.
By Lemma~\ref{E-def.thm}, $\E$ is isomorphic to $\uGamma_Y(\Bbb I_X)$ in
$D(X)$.
On the other hand, $\uGamma_Y(\Bbb I_X)$ is quasi-isomorphic to 
$\Bbb J=\Cone(\Bbb I_X\rightarrow f_*f^*\Bbb I_X)[-1]$, where $f:X\setminus Y
\rightarrow X$ is the inclusion.
As $f_*f^*$ has an exact left adjoint $f_!f^*$ 
(see (\ref{preserving-injectives.par})), $\Bbb J$ is a bounded
injective resolution of $\E$.
\qed

\begin{lemma}\label{vanishing.thm}
$\uExt_{\O_X}^i(\M,\E)=0$ for $i>0$ and $\M\in\Coh(G,X)$.
\end{lemma}

\proof 
$\uExt_{\O_X}^i(\M,\E)_\eta\cong \Ext^i_{\O_{X,\eta}}(\M_\eta,\E_\eta)$.
As $\E_\eta$ is injective, we are done.
\qed

\begin{corollary}
$\Bbb D:=\uHom_{\O_X}(?,\E)$ is an exact functor on $\Coh(G,X)$.
\end{corollary}

\begin{lemma}\label{keroron.thm}
For $\M\in\Qch(G,X)$, the following are equivalent:
\begin{description}
\item[(i)] $\M$ is of finite length;
\item[(ii)] $\M\in\Coh(G,X)$, and $\I^n\M=0$ for some $n$.
\item[(iii)] $\M_\eta$ is an $\O_{X,\eta}$-module of finite length;
\end{description}
\end{lemma}

\proof {\bf (i)$\Rightarrow$(ii)} 
As $\M$ is of finite length, it is a noetherian object.
Hence it is coherent by \cite[Lemma~12.8]{ETI}.
As $\M$ is also an artinian object, $\I^n\M=\I^{n+1}\M$ for sufficiently
large $n$.
By Corollary~\ref{G-NAK}, $\I^n\M=0$.

{\bf (ii)$\Rightarrow$(iii)}
As $\M$ is coherent, $\M_\eta$ is a finitely generated $\O_{X,\eta}$-module.
Since $\I_\eta^n\M_\eta=0$, the support of $\M_\eta$ is one point,
and hence $\M_\eta$ is a module of finite length.

{\bf (iii)$\Rightarrow$(i)} This is because $(?)_\eta$ is faithfully exact.
\qed

\paragraph We denote by $\Cal F$ the full subcategory of those objects 
$\M\in\Qch(G,X)$ such that the equivalent conditions in the lemma
are satisfied.

\begin{theorem}[Matlis duality]\label{matlis.thm}
Set $\Bbb D:=\uHom_{\O_X}(?,\E)$.
Then
\begin{description}
\item[(i)] $\Bbb D$ is an exact functor from $\Cal F$ to itself.
\item[(ii)] $\Bbb D^2\cong \Id$ as functors on $\Cal F$.
In particular, $\Bbb D:\Cal F\rightarrow \Cal F$ is an anti-equivalence.
\end{description}
\end{theorem}

\proof {\bf (i)}
If $\M\in\Cal F$, then $\Bbb D(\M)=\uHom_{\O_X}(\M,\E)$ is 
in $\Qch(G,X)$, and
$\uHom_{\O_X}(\M,\E)_\eta=\Hom_{\O_{X,\eta}}
(\M_\eta,\E_\eta)$ is of finite length,
because this module is the Matlis dual of the module $\M_\eta$, which is
of finite length.
So the condition {\bf (iii)} in Lemma~\ref{keroron.thm} is satisfied, and
hence $\Bbb D(\M)\in\Cal F$.
The exactness of $\Bbb D$ is already checked.

{\bf (ii)} Let $\fD:\M\rightarrow \Bbb D\Bbb D\M=
\uHom_{\O_X}(\uHom_{\O_X}(\M,\E),\E)$ be the canonical map, see (\ref{fD.par}).
Note that by Lemma~\ref{fD-H.thm} and Lemma~\ref{fD-P.thm}, 
applying $(?)_\eta$ to this map, we get the duality map
$\fD:\M_\eta\rightarrow \Hom_{\O_{X,\eta}}(\Hom_{\O_{X,\eta}}(
\M_\eta,\E_\eta),\E_\eta)$, which is an isomorphism, since $\E_\eta$ is
the injective hull of the residue field $\kappa(\eta)$.
Since $(?)_\eta$ is faithful, $\fD: \M\rightarrow\Bbb D\Bbb D\M$
is an isomorphism.
\qed

\begin{theorem}[Local duality]\label{local-duality.thm}
For $\Bbb F\in D_{\Coh}(G,X)$, 
the composite
\begin{multline*}
\fd:R\uGamma_Y \Bbb F
\xrightarrow{\fD}
R\uGamma_Y R\uHom_{\O_X}(R\uHom_{\O_X}(\Bbb F,\Bbb I_X),\Bbb I_X)
\xrightarrow{\fH}\\
R\uHom_{\O_X}(R\uHom_{\O_X}(\Bbb F,\Bbb I_X),R\uGamma_Y\Bbb I_X)
\cong
R\uHom_{\O_X}(R\uHom_{\O_X}(\Bbb F,\Bbb I_X),\E)
\end{multline*}
is an isomorphism.
It induces an isomorphism
\[
\uH^i_Y(\Bbb F)\cong \uHom_{\O_X}(\uExt^{-i}_{\O_X}(\Bbb F,\Bbb I_X),\E)
\]
for each $i\in\Bbb Z$.
\end{theorem}

\proof $\fD$ in the composition is an isomorphism by
\cite[(31.9)]{ETI}.
$\fH$ is an isomorphism by Theorem~\ref{fH-isom-derived.thm}, {\bf (i)}.
Thus $\fd$ is an isomorphism.

To prove the second assertion,
it suffices to show that
\[
\uExt^i_{\O_X}(\Bbb G,\E)\cong \uHom_{\O_X}(H^{-i}(\Bbb G),\E),
\]
where $\Bbb G=R\uHom_{\O_X}(\Bbb F,\Bbb I_X)$.
Note that $\Bbb G\in D_{\Coh}(X)$ by \cite[(31.9)]{ETI}.
Let $\Bbb J$ be a bounded injective resolution of $\E$ (it does exist, 
see Lemma~\ref{fid.thm}).
Consider the spectral sequence
\[
E^{p,q}_2=H^p(\uHom_{\O_X}(H^{-q}(\Bbb G),\Bbb J))\Rightarrow \uExt^{p+q}_
{\O_X}(\Bbb G,\E).
\]
By Lemma~\ref{vanishing.thm}, $E^{p,q}_2=0$ for $p\neq 0$, and the 
spectral sequence collapses, and we get the desired assertion.
\qed

\begin{lemma}
Let $\Bbb F\in D_{\Coh}(X)$.
Then the diagram
\[
\xymatrix{
(?)_\eta(?)_0 R\uGamma_Y(\Bbb F)
\ar[r]^-\fd \ar[dd]^{\bar \delta \hat\gamma^{-1}} &
(?)_\eta(?)_0R\uHom_{\O_X}(R\uHom_{\O_X}(\Bbb F,\Bbb I_X),\E)
\ar[d]^{PH} \\
 &
R\Hom_{\O_{X,\eta}}((?)_\eta(?)_0 R\uHom_{\O_X}(\Bbb F,\Bbb I_X),\E_\eta)
\ar[d]^{P^{-1}H^{-1}}
\\
R\Gamma_{\fm_\eta}(\Bbb F_\eta)
\ar[r]^-\fd &
R\Hom_{\O_{X,\eta}}(R\Hom_{\O_{X,\eta}}(\Bbb F_\eta,\Bbb I_{X,\eta}),\E_\eta)
}
\]
is commutative \(see for the definition of $\hat\gamma$ and $\bar\delta$, 
see {\rm \cite[section~4]{HO}} and {\rm
\cite[(6.1)]{HO}}, respectively\).
\end{lemma}

\proof Note that $H^{-1}$ in the diagram exists by
\cite[(13.9)]{ETI}.
The $P^{-1}$ exists by Lemma~\ref{derived-P.thm}, {\bf (iii)}.
The commutativity of the diagram follows from
Lemma~\ref{fD-H.thm} and Lemma~\ref{fD-P.thm} immediately.
\qed

\section{An example of graded rings}

\paragraph
Let $(R,\fm)$ be a noetherian local ring with a normalized dualizing complex
$\Bbb I_R$.
Set $S=\Spec R$.
Let $H$ be a flat $R$-group scheme of finite type, and $G=\Bbb G_m\times H$.
Let $A$ be a $G$-algebra.
So $A$ is $\Bbb Z$-graded and each homogeneous component is an $H$-submodule
of $A$.
Assume that $A=\bigoplus_{i\geq 0}A_i$ is $\Bbb N$-graded and
$A_0=R$.
Let $\pi:X\rightarrow S$ be the canonical map, where $X:=\Spec A$.
Set $\Bbb I_X:=\pi^!\Bbb I_R$.

\begin{lemma}
Under the notation as above, $X$ is $G$-local, and 
$\Bbb I_X$ is $G$-normalized.
\end{lemma}

\proof Let $I$ be a proper $G$-ideal of $A$.
Then $I$ is a homogeneous ideal, and is contained in the unique graded 
maximal ideal $\fM:=\fm+A_+$, where $A_+=\bigoplus_{i>0}A_i$.
Clearly, $\fM$ is a $G$-ideal, and hence is the unique $G$-maximal 
$G$-ideal.
So $X$ is $G$-local.

Let $\varphi: S\rightarrow X$ be the closed immersion induced by
$A\rightarrow A/A_+=R$.
Let $\psi:Y\rightarrow S$ be the closed immersion induced by
$R\rightarrow R/\fm\cong A/\fM$, where $Y=\Spec A/\fM$.
Then since $\pi\varphi=\id_S$, 
\[
\Bbb I_Y=(\varphi\psi)^!(\Bbb I_X)=\psi^!\varphi^!\pi^!\Bbb I_R
=\psi^!\Bbb I_R.
\]
So $H^i(\Bbb I_Y)\cong \Ext_R^i(R/\fm,\Bbb I_R)$, whose Matlis dual is
$H^i_\fm(R/\fm)$.
This is nonzero if and only if $i=0$.
Thus $\Bbb I_X$ is $G$-normalized.
\qed

\paragraph For a finite $R$-module $V$, set $V^\dagger:=\Hom_R(V,E_R)$, 
where $E_R$ is the injective hull of the residue field $R/\fm$ of $R$.
For an $A$-finite $G$-module $M$, set 
$M^\vee=\indlim\Hom_R(M/\fM^n M,E_R)$.
As each $M/\fM^n M$ is an $R$-finite $(G,A)$-module, each
$\Hom_R(M/\fM^n M,E_R)$ is a $(G,A)$-module, and hence $M^\vee$ is
also a $(G,A)$-module.
It is easy to see that $M^\vee\cong \Hom_A(M,A^\vee)$.
Note that the degree $i$ component of $M^\vee $ is $M_{-i}^\dagger$.
That is, $M^\vee=\bigoplus_{i\in\Bbb Z}M_{-i}^\dagger$.

\begin{lemma} $A^\vee$ is isomorphic to $E_A:= \Gamma(X,\E)$ 
as a $(G,A)$-module.
\end{lemma}

\proof We may assume that $\Bbb I_R$ is the normalized fundamental 
dualizing complex.
We have
\begin{multline}\label{indlim.eq}
\E=\uH^0_{Y}(\Bbb I_X)=\indlim \uExt_{\O_X}^0(\O_X/\I^n,\Bbb I_X)\\
=\indlim H^0((\psi_n)_*\psi_n^!\pi^!\Bbb I_S)
=
\indlim \Ext_R^0(A/\fM^n,\Bbb I_R)\,\widetilde{\;},
\end{multline}
where $\psi_n:\Spec A/\fM^n\rightarrow \Spec A$ is the canonical
closed immersion,
and $(?)\,\widetilde{\;}$ denotes the quasi-coherent sheaf associated 
to a module.
On the other hand, $A/\fM^n$ has finite length as an $R$-module, 
so
\begin{multline*}
\Ext_R^0(A/\fM^n,\Bbb I_R)
\cong H^0(\Hom_R(A/\fM^n,\Bbb I_R))
\cong H^0(\Hom_R(A/\fM^n,\Gamma_{\fm}\Bbb I_R))\\
\cong H^0(\Hom_R(A/\fM^n,E_R))=\Hom_R(A/\fM^n,E_R).
\end{multline*}

We prove that 
the map $\Hom_R(A/\fM^m,E_R)\rightarrow \Hom_R(A/\fM^n,E_R)$
in the inductive system is induced by the projection $A/\fM^n\rightarrow
A/\fM^m$ for $n\geq m$.
Note that $\uExt_{\O_X}^0(\O_X/\I^m,\Bbb I_X)
\rightarrow\uExt_{\O_X}^0(\O_X/\I^n,\Bbb I_X)$ 
in (\ref{indlim.eq}) is induced by the projection.
So by the description of the twisted inverse for finite morphisms
\cite[(27.7)]{ETI},  the map $(\psi_m)_*\psi_m^!\rightarrow (\psi_n)_*
\psi_n^!$ is induced by the counit map.
That is, 
the map is the composite
\[
(\psi_m)_*\psi_m^!\cong
(\psi_n)_*(\psi_{n,m})_*\psi_{n,m}^!\psi_n^!
\xrightarrow{\varepsilon}
(\psi_n)_*\psi_n^!,
\]
where $\psi_{n,m}:\Spec A/\fM^m\rightarrow \Spec A/\fM^n$ is the map
induced by the projection.
So again by \cite[(27.7)]{ETI}, the map
$\Ext_R^0(A/\fM^m,\Bbb I_R) \rightarrow \Ext_R^0(A/\fM^n,\Bbb I_R)$ in
(\ref{indlim.eq}) is also induced by the projection, and we are done.

Hence
\[
E_A=\indlim \Hom_R(A/\fM^n,E_R)=A^\vee.
\qed
\]

\begin{corollary}\label{graded-local-duality.thm}
Assume that $A$ is Cohen--Macaulay and $\dim A=d$.
Set $\Gamma(X,H^{-d}(\Bbb I_X))$ to be $\omega_A$.
For a $A$-finite $(G,A)$-module $M$, the canonical map
\[
\fd: H^i_{\fM}(M)\rightarrow \Ext^{d-i}_A(M,\omega_A)^\vee
\]
is an isomorphism of $(G,A)$-modules.
That is, this isomorphism preserves grading and $H$-action.
\end{corollary}

\begin{remark}\label{equivariant-Matlis-graded.thm}
Assume that $R=k$ is a field.
Let $\Cal G$ be the full subcategory of $(G,A)$-modules
consisting of $M$ such that $M_i$ is finite dimensional for every $i$.
Then we define $M^\vee=\bigoplus_{i\in\Bbb Z}M_{-i}^\dagger$ for $M\in 
\Cal G$, where $M_{-i}^\dagger=\Hom_k(M_{-i},k)$.
We have an isomorphism $\Phi: \sHom_A(M,A^\vee)\rightarrow M^\vee$.
See for the notation $\sHom_A$, \cite{BH}.
Note that
\begin{multline*}
\Phi_n:\sHom_A(M,A^\vee)_n=\sHom_{A}(M(-n),A^\vee)_0\\
\rightarrow \Hom_k(M_{-n},A^\vee_0)=\Hom_k(M_{-n},k)
\end{multline*}
is given by the restriction.
It is easy to see that $(?)^\vee$ is an anti-equivalence from $\Cal G$
to itself.
This also gives an anti-equivalence between the category of noetherian
$(G,A)$-modules to that of artinian $(G,A)$-modules.
This is not contained in Theorem~\ref{matlis.thm}, which treats only
objects of finite length.
\end{remark}

\begin{example}\label{non-injective.thm}
Let $k$ be an algebraically closed field of characteristic two,
and we set $R=k$ and $S=\Spec R$.
Let $V=k^2$, and $H=\GL(V)$.
Let $A=\Sym V$, and $X=V^*=\Spec A$.
Then $A_2^*$ is not injective as a $G$-module.
So $E_A=\bigoplus_{i\geq 0}A_i^*$ is not injective as a $G$-module either.
So $E_A$ is not injective as a $(G,A)$-module either by 
\cite[Corollary~II.1.1.9]{Hashimoto}.
In particular, $E_A$ is not the injective hull of $A/\fM$ as a $(G,A)$-module.
\end{example}

\end{document}